\documentclass[centertags,12pt]{amsart}
\usepackage{latexsym}
\usepackage{amssymb}

\numberwithin{equation}{section}

\makeatletter
\renewcommand{\subsection}{\@startsection
{subsection}{2}{0mm}{\baselineskip}{-0.25cm}{\normalfont\normalsize\rm}}
\makeatother

\newtheorem{theorem}{Theorem}[section]
\newtheorem{proposition}[theorem]{Proposition}
\newtheorem{lemma}[theorem]{Lemma}
\newtheorem{corollary}[theorem]{Corollary}
{\theoremstyle{definition}

\newtheorem*{conjecture*}{Conjecture}}
\theoremstyle{remark}
\newtheorem{remark}[theorem]{Remark}
\newtheorem{claim}{Claim}

\def\F{\mathbf F}
\def\N{\mathbf N}
\def\P{\mathbf P}

\def\cD{\mathcal D}
\def\cH{\mathcal H}

\def\cK{\mathcal K}

\def\cQ{\mathcal Q}

\def\cX{\mathcal X}

\def\fl{\mathbf{F_\ell}}
\def\fq{\mathbf F_{q^2}}
\def\dim{{\rm dim}}

\def\deg{{\rm deg}}

\def\frx{{\mathbf{Fr}}_{\mathcal X}}

\def\supp{{\rm Supp}}
\begin{document}

\author[G.~Korchm\'aros]{G\'abor Korchm\'aros}
\author[F.~Torres]{Fernando Torres}
\thanks{This research was carried out
within the project ``Progetto e Realizzazione di un Criptosistema per
Telecomunicazioni", POP FESR 1994/99 - II Triennio. The second author
was partially support by Cnpq-Brazil, Proc. 300681/97-6. We are indebted
to E. Ballico for sending us his paper \cite{ballico}}

\title[The genus of a maximal curve]{On the genus of a maximal curve}

   \address{Dipartimento di Matematica, Universit\'a della
Basilicata, via N. Sauro 85,\\
85100 Potenza, Italy}
   \email{korchmaros@unibas.it}
   \address{IMECC-UNICAMP, Cx. P. 6065, Campinas-13083-970-SP, Brazil}
   \email{ftorres@ime.unicamp.br}

    \begin{abstract} Previous results on genera $g$ of $\fq$-maximal curves
are improved:
  \begin{enumerate}
\item[\rm(1)] $\text{Either} \ g\leq \lfloor (q^2-q+4)/6\rfloor\,,
\ \text{or} \ g=\lfloor(q-1)^2/4\rfloor\,, \ \text{or} \
g=q(q-1)/2\,$.

\item[\rm(2)] The hypothesis on the existence of a particular
Weierstrass point in \cite{at} is proved.

\item[\rm(3)] For $q\equiv 1\pmod{3}$, $q\ge 13$, no $\fq$-maximal curve
of genus $(q-1)(q-2)/3$ exists.

\item[\rm(4)] For $q\equiv 2\pmod{3}$, $q\ge 11$, the non-singular
$\fq$-model of the plane curve of equation $y^q+y=x^{(q+1)/3}$ is the
unique $\fq$-maximal curve of genus $g=(q-1)(q-2)/6$.

\item[\rm(5)] Assume $\dim(\cD_\cX)=5$, and ${\rm char}(\fq)\geq 5$. For
$q\equiv 1\pmod{4}$, $q\geq 17$, the Fermat curve of equation
$x^{(q+1)/2}+y^{(q+1)/2}+1=0$ is the unique $\fq$-maximal curve of genus
$g=(q-1)(q-3)/8$. For $q\equiv 3\pmod{4}$, $q\ge 19$, there are exactly
two $\fq$-maximal curves of genus $g=(q-1)(q-3)/8$, namely the above
Fermat curve and the non-singular $\fq$-model of the plane curve of
equation $y^q+y=x^{(q+1)/4}$. 
   \end{enumerate}
The above results provide some new evidences on maximal curves in
connection with Castelnuovo's bound and Halphen's theorem,
especially with extremal curves; see for instance the conjecture
stated in Introduction.
    \end{abstract}

\maketitle

\section{Introduction}\label{s1}

An {\em $\fq$-maximal curve} $\cX$ of genus $g$ is defined to be a
projective, geometrically irreducible, non-singular algebraic curve
defined over $\fq$ such that the number of its $\fq$-rational points
attains the Hasse-Weil upper bound, namely
  $$
\#\cX(\fq)=q^2+1+2qg\, .
  $$
$\fq$-maximal curves especially those with large genus are currently
investigated also in connection with coding theory and cryptography based
on Goppa's method \cite[Ch. 4, Sect. 7]{goppa}. It is well known that
$g\leq q(q-1)/2$, see \cite{ihara}, and that $g$ reaches this upper limit
if and only if $\cX$ is $\fq$-isomorphic to the Hermitian curve, see
\cite{r-sti}. In \cite{ft1} it is proven that
   \begin{equation}\label{ft1}
\text{either}\qquad g \leq \lfloor(q-1)^2/4\rfloor\, ,
\qquad\text{or}\qquad g=q(q-1)/2\, .
    \end{equation}
For $q$ odd, $g=(q-1)^2/4$ occurs if and only if $\cX$ is $\fq$-isomorphic
to the non-singular model of the plane curve of equation
$y^q+y=x^{(q+1)/2}$, see \cite[Thm. 3.1]{fgt}. For $q$ even, a similar
result is obtained in \cite{at} under an extra-condition that $\cX$ has a
particular Weierstrass point: $g=\lfloor (q-1)^2/4\rfloor=q(q-2)/4$ if and
only if $\cX$ is $\fq$-isomorphic to the non-singular model of the plane
curve of equation $y^{q/2}+\ldots+y^2+y=x^{q+1}$. These results together
with some evidences coming from \cite{ckt1}, \cite{ckt2}, \cite{g-sti-x}
make it plausible that only few $\fq$-maximal curves can have genus close
to the upper limit. As a matter of fact, in the range
   $$
\lfloor (q-1)(q-3)/8\rfloor\leq g<\lfloor (q-1)^2/4\rfloor\, ,
   $$
only twelve examples up to $\fq$-isomorphisms are known to exist and the
spectrum of their genera is listed below:
   \begin{enumerate}
\item[\rm(I)] $g=\lfloor (q^2-q+4)/6\rfloor$ for $q\equiv 0,1,2\pmod{3}$,
see Remark \ref{rem3.2};

\item[\rm(II)] $g=(q^2-q-2)/6$ for $q\equiv 2\pmod{3}$, see \cite[Thm.
6.2]{ckt1} or \cite[Thm. 5.1]{g-sti-x};

\item[\rm(III)] $g=\lfloor ((q-1)(q-2)/6\rfloor$ for $q\equiv
0,2\pmod{3}$, see the case $N=4$ in (\ref{eq2.331}), and Sect. \ref{s4.1};

\item[\rm(IV)] $g=\lfloor (q^2-2q+5)/8\rfloor$ for $q\equiv
0,1,3\pmod{4}$, see Remark \ref{rem4.21};

\item[\rm(V)] $g=\lfloor (q-1)(q-3)/8\rfloor$ for $q\equiv 0,1,3\pmod{4}$,
see the case $N=5$ in (\ref{eq2.331}), and Sect. \ref{s4.2}.
   \end{enumerate}
Theorem \ref{thm3.1} in this paper together with (\ref{ft1}) prove
the following result, see Corollary \ref{cor3.1}:
   \begin{equation}\label{kt}
\text{either}\quad g\leq \lfloor (q^2-q+4)/6\rfloor\, ,\quad
\text{or}\quad g=\lfloor(q-1)^2/4\rfloor\, ,\quad\text{or}\quad
g=q(q-1)/2\, .
   \end{equation}
This result is the best possible since the upper bound in
(\ref{kt}) cannot be improved as it is attained by the curves cited
in (I) for every $q$. In other words the third largest genus of an
$\fq$-maximal curve equals $g=\lfloor (q^2-q+4)/6\rfloor$
independently of $q$; by contrast, the fourth largest genus might
heavily depend on $q$. The above examples also show that the gap
between the first and second as well as the second and third
largest genus is approximately constant times $q^2$, while the gap
between the third and forth is only $1$ for $q\equiv 2\pmod{3}$,
and at most constant times $q$ for $q\equiv 0\pmod{3}$.

The essential idea of the proof of Theorem \ref{thm3.1} is to show
that every $\fq$-maximal curve of genus $g>\lfloor
(q^2-q+4)/6\rfloor$ has a non-singular model $\cX$ over $\fq$
embedded in $\P^3(\bar\fq)$ such that $\cX$ has degree $q+1$ and
lies on an $\fq$-rational quadratic cone $\cQ$ whose vertex $V$
belongs to $\cX$. This idea will be worked out using the ``natural
embedding theorem", see \cite[Thm. 2.5]{kt1}, together with
Weierstrass point theory, Castelnuovo's genus bound, Halphen's
theorem and some other tools. Actually, for $q$ even the vertex $V$
is a particular $\fq$-rational Weierstrass point of $\cX$, since
the order-sequence of $\cX$ at $V$ (i.e. the possible intersection
numbers of $\cX$ with hyperplanes at $V$) turns out to be
$(0,1,(q+2)/2,q+1)$. Similarly for $q$ odd, an $\fq$-rational
Weierstrass point with order-sequence $(0,1,(q+1)/2,q+1)$ is shown
to exist. An $\fq$-rational point of $\cX$ with such a particular
order-sequence forces $\cX$ to have genus $\lfloor(q-1)^2/4\rfloor$
as noticed in \cite[Remark 2.6(1)]{kt1}. Then, the already quoted
characterization theorems from \cite{fgt}, and \cite{at} will be
applied to complete the proof of Theorem \ref{thm3.1}. It should be
noted that Theorem \ref{thm3.1} improves both \cite[Prop. 2.5]{ft2}
and the main result in \cite{at}.

Curves with genera as in (III) and (V) above turn out be {\em
extremal} in $\P^4(\bar\fq)$ and in $\P^5(\bar\fq)$ respectively,
as such genera are Castelnuovo's numbers $c_0(q+1,r)$, $r=4,5$, see
(\ref{castelnuovo}). Extremal curves in zero characteristic have
been widely investigated, see, for instance, \cite{acgh} and the
references therein. Several relevant properties of extremal curves
are known to hold true in positive characteristic thanks to
Rathmann's work \cite{rathmann} (see also
\cite{ballico-cossidente}). For the present purpose, the key result
on extremal curves is Lemma \ref{lemma2.12} stated in Sect.
\ref{s2.1}. Indeed, this lemma together with other results will
give both the non-existence of $\fq$-maximal curves of genus
$(q-1)(q-2)/3$ for $q\equiv 1\pmod{3}$, and a characterization of a
$\fq$-maximal curve with such a genus for $q\equiv 2\pmod{3}$,
$q\ge 11$; see Theorem \ref{thm4.11}. Under two additional
hypotheses, namely the curve is naturally embedded in
$\P^5(\bar\fq)$ and $char (\fq) \geq 5$, the aforementioned key
result will also be an essential ingredient in characterizing
$\fq$-maximal curves of genus $(q-1)(q-3)/8$ for $q\equiv
1,3\pmod{4}$, $q\ge 11$; see Theorem \ref{thm4.21}. This theorem is
related to a previous characterization of plane $\fq$-maximal
curves of degree $(q+1)/2$ stated in \cite{chkt}. Also, in view of
the results in Sect. \ref{s4.1} and \cite[Proof of Prop. 2.4]{at},
it seems plausible that any two $\fq$-maximal curves of genus
$q(q-3)/6$ for $q\equiv 0\pmod{3}$ are $\fq$-isomorphic. On the
contrary, due to the examples in \cite[Sect. 5]{abdon}, no similar
result for curves of genus $q(q-4)/8$, $q\equiv 0\pmod{4}$ can
hold. For a further interesting question related to these matters,
see Remark \ref{rem2.31}.

The genera in (I) and (IV) above coincide with Halphen's number
$c_1(q+1,r)$, $r=3,4$, see (\ref{halphen}). Extensions of results
around Halphen's theorem from zero characteristic to positive
characteristic are also possible again by Rathmann's work
\cite{rathmann} and Ballico's paper \cite{ballico}. Unfortunately,
we do not have so far a classification theorem for $\fq$-maximal
curves with such genera. What we currently know in this direction
is that extremal curves lie on special surfaces such as scrolls,
see e.g. \cite[Ch. III, Thm. 2.5]{acgh}, and that curves with
enough large degree and genus equal to Halphen's number are
Cohen-Macaulay curves lying on Castelnuovo surfaces, see the main
theorem in \cite{chi-ci-ge}. These facts together with the general
form of the above mentioned ``natural embedding theorem'' stating
that every $\fq$-maximal curve is naturally embedded in a
high-dimensional projective space over $\fq$ as a curve of degree
$q+1$ contained in a Hermitian variety of degree $q+1$, see
\cite[Thm. 3.4]{kt1}, seem to be a good starting point of a
classification project for such $\fq$-maximal curves.

Finally, we stress that (\ref{kt}) provides evidence for the
following conjecture.
   \begin{conjecture*} With notation as in (\ref{castelnuovo}) and
(\ref{halphen}), there is no $\fq$-maximal curve of
genus $g$ such that
   $$
c_1(q+1,r)<g<c_0(q+1,r)\, .
   $$
   \end{conjecture*}

\section{Background}\label{s2}

Our purpose in this section is to recall some results concerning upper
bounds on the genus of curves, Weierstrass Point Theory and Frobenius
orders as well as some results on maximal curves.

{\bf Convention.} The word {\em curve} will mean a projective,
geometrically irreducible, non-singular algebraic curve.

\subsection{Castelnuovo's genus bound and Halphen's theorem}\label{s2.1}

Throughout this sub-section, $\cX$ denotes a curve defined over an
algebraically closed field $\F$. Let $\cD$ be an $r$-dimensional, $r\ge
2$, base-point-free linear series of degree $d$ defined on $\cX$; $\cD$ is
assumed to be simple, that is $\cX$ is birational to $\pi(\cX)$, where
$\pi$ denotes a morphism associated to $\cD$. Castelnuovo showed that the
genus $g$ of $\cX$ is upper bounded by a function depending on $r$ and
$d$. More precisely, let $\epsilon$ be the unique integer with
$0\le\epsilon\le r-2$ and $d-1\equiv \epsilon \pmod{(r-1)}$, and define
Castelnuovo's number $c_0(d,r)$ by
  \begin{equation}\label{castelnuovo}
c_0(d,r):=\frac{d-1-\epsilon}{2(r-1)}(d-r+\epsilon)\, .
  \end{equation}

  \begin{lemma} {\rm (Castelnuovo's genus bound for curves in projective
spaces, \cite{castelnuovo}, \cite[p. 116]{acgh}, \cite[IV, Thm.
6.4]{hartshorne}, \cite[Thm. 3.3]{accola}, \cite[Cor.
2.8]{rathmann})}\label{lemma2.11}
  $$
g\le c_0(d,r,\epsilon)\, .
  $$
  \end{lemma}
  \begin{remark}\label{rem2.11}
  $$
c_0(d,r,\epsilon)\le
\begin{cases}
(d-1-(r-1)/2)^2/2(r-1) & \text{for $r$ odd,}\\
(d-1-(r-1)/2)^2-1/4)/2(r-1) & \text{for $r$ even.}
\end{cases}
 $$
  \end{remark}
Curves with genus equal to Castelnuovo's number have several remarkable
properties; see e.g. \cite{accola}, \cite[Ch. 3]{eisenbud-harris},
\cite[Ch. 3, Sect. 2]{acgh}. We will use the following one, which
is in fact implicitly contained in the proof of Castelnuovo's genus bound
taking into account the Riemann-Roch theorem; see e.g. \cite[p. 361 and Lemma
3.5]{accola}.
  \begin{lemma}\label{lemma2.12} Assume $g=c_0(d,r)$, and define
$\epsilon'$ by $d=m(r-1)+\epsilon'$ with $\epsilon'\in\{2,\ldots, r\}.$ If
$m\ge 2.$, then:
  \begin{enumerate}
  \item[\rm(1)] the dimension of the linear series $2\cD$ is $3r-1;$
  \item[\rm(2)] there exists a base-point-free $(\epsilon'-2)$-dimensional
complete linear series $\cD'$ of degree $(\epsilon'-2)(m+1)$ such that
$(m-1)\cD+\cD'$ is the canonical linear series$.$
\end{enumerate}
  \end{lemma}
The following theorem going back to Halphen improves Castelnuovo's genus
bound for certain curves in $\P^3(\F)$.
  \begin{lemma} {\rm (Halphen's theorem, \cite[Thm. 3.1]{gruson-peskine},
\cite[Thm.  3.13]{eisenbud-harris}, \cite{ballico})}\label{lemma2.13}
Assume $d\ge 7$, and $d=17$ or $d\ge 25$ when $char(\F)>0$. If $\cX$ is
embedded in $\P^3(\F)$, then $\cX$ lies on a quadric surface provided that
   $$
g> c_1(d,3):=\lfloor (d^2-3d+6)/6\rfloor\, .
   $$
   \end{lemma}
   \begin{remark}\label{rem2.12} For a historical account of
   Halphen's theorem, see \cite[p. 349]{hartshorne} or Introduction in
\cite{gruson-peskine} and \cite{harris}. The proof in characteristic $0$
due to Eisenbud and Harris \cite[Thm. 3.13]{eisenbud-harris} depends on
the Uniform Position Principle applied to the generic hyperplane section
of $\cX$, and it still works verbatim in positive characteristic.
   \end{remark}
Halphen's theorem extends to certain curves in $\P^r(\F)$ for $r\ge
4$, and it turns out to be very useful when one looks for a bound
$c_\alpha(d,r)$ for the genus of a curve of degree $d$ in
$\P^r(\F)$ not lying on any irreducible surface of degree less than
$r+\alpha-1$. For our purpose, the smallest case $\alpha=1$ is
needed:
   \begin{lemma}\label{lemma2.14} {\rm (\cite[Thm. 3.22]{eisenbud-harris},
\cite[Cor. 2.8]{rathmann})} Suppose that $\cX$ is a curve in $\P^r(\F)$ of
degree $d$ and genus $g$. Assume
  $$
d\ge \begin{cases}
36r   &  \text{if $r\le 6$,}\\
288   &  \text{if $r=7$,}\\
2^{r+1} & \text{if $r\ge 8$.}\end{cases}
  $$
Then $\cX$ lies on a surface of degree less than
or equal to $r-1$ provided that
   \begin{equation}\label{halphen}
g>c_1(d,r):=\frac{d-1-\epsilon_1}{2r}(d-r+\epsilon_1+1)+
\begin{cases}
0 & \text{if $\epsilon_1\le r-2$}\\
1 & \text{if $\epsilon_1=r-1$}\end{cases}\, ,
   \end{equation}
where $\epsilon_1$ is the unique integer such that $0\le \epsilon_1\le
r-1$ and $d-1\equiv \epsilon_1\pmod{r}.$
   \end{lemma}
Notice that (\ref{halphen}) for $r=3$ coincides with the formula in
Lemma \ref{lemma2.13}. A full account of results related to
Halphen's theorem is found in the already mentioned
\cite{eisenbud-harris}, \cite{gruson-peskine}, as well as in
\cite{chi-ci-ge} and \cite{chi-ci}.

\subsection{Weierstrass Point Theory and Frobenius orders}\label{s2.2}

Our reference in this sub-section is St\"ohr-Voloch's paper \cite{sv}.
Let $\cX$ be a curve defined over an algebraically closed field $\F$ of
characteristic $p$, $g$ its genus, and $\cD$ an $r$-dimensional, $r\ge 1$,
simple base-point-free linear series of degree $d$ defined on $\cX$. The
{\em $(\cD,P)$-order sequence} of $P\in\cX$ is the strictly increasing
sequence $j_0(P)=0<j_1(P)<\ldots<j_r(P)$ enumerating the set
$\{v_P(D):D\in\cD\}$, with $v_P(D)$ being the weight of the divisor $D$ at
$P$, see \cite[p. 3]{sv}. If $\pi$ is a morphism associated to $\cD$, then
   $$
\cD=\{\pi^*(H):\, \text{$H$ hyperplane in $\P^r(\F)$}\}\, ,
   $$
and the $(\cD,P)$-order sequence consists of all possible intersection
numbers of $\cX$ with hyperplanes at $P$ in the usual order whenever
$\cX\subseteq \P^r(\F)$. Furthermore, the $(\cD,P)$-order sequence is the
same for all but finitely many points \cite[pp. 4-6]{sv}. and each of the
exceptional points is called a {\em $\cD$-Weierstrass point} of $\cX$.
According to \cite[p. 6]{sv}, there exists a divisor $R=R_{\cD}$ on $\cX$,
the so-called {\em ramification divisor}, with support consisting of all
$\cD$-Weierstrass points of $\cX$ and degree
  \begin{equation}\label{eq2.21}
\deg(R)=\sum_{i=0}^{r}\epsilon_i(2g-2)+(r+1)d\, ,
  \end{equation}
where $\epsilon_0=0<\epsilon_1=1<\ldots<\epsilon_r$ is the {\em
$\cD$-order sequence} of $\cX$, that is the $(\cD,P)$-order
sequence at a general (i.e. a non $\cD$-Weierstrass) point
$P\in\cX$. It should be noted that the well known {\em Weierstrass
points of $\cX$} appear in this context as the Weierstrass points
of the canonical linear series on $\cX$ in which case
   $$
H(P):=\N_0\setminus\{j_i(P)+1:i=0,\ldots,g-1\}
   $$
is a numerical semigroup whose elements are called {\em Weierstrass
non-gaps} at $P$. The strictly increasing sequence enumerating
$H(P)$ is usually denoted by $(m_i(P):i=0,1,\ldots)$.

A general rule to compute the $(\cD,P)$-orders and $v_P(R)$ is given by
the following lemma.
   \begin{lemma} {\rm (\cite[p. 5, Thm. 1.5]{sv})}\label{lemma2.21}
\begin{enumerate}
\item[\rm(1)] $j_i(P)\ge \epsilon_i$ for each $P$ and each $i;$
\item[\rm(2)] $v_P(R)\ge \sum_{i=0}^{r}(j_i(P)-\epsilon_i),$ and equality
holds if and only if ${\rm det}(\binom{j_i(P)}{\epsilon_k})\not\equiv
0\pmod{p}.$
\end{enumerate}
    \end{lemma}
To every point $P\in\cX$ there is attached the flag of osculating
subspaces of $\P^r(\F)$ relative to a morphism $\pi$ associated
to $\cD$. For each $i$, $0\leq i\leq r-1$, the {\em $i$th
osculating space} $L_i(P)$ of $\cX$ at $P$ (with respect to $\pi$)
is the $i$-dimensional subspace in $\P^r(\F)$ defined as the
intersection of all hyperplanes $H$ in $\P^r(\F)$ satisfying
$v_P(\pi^*(H))\ge j_{i+1}(P)$. Clearly, $L_0(P)=\{\pi(P)\}\subseteq
L_1(P)\subseteq\ldots\subseteq L_{r-1}(P)$. Also, $L_i(P)$ is
uniquely determined by $\cD$ up to projective equivalence because
any two morphisms associated to $\cD$ are projectively equivalent.
We will refer to $L_1(P)$ and $L_{r-1}(P)$ as the {\em tangent
line} and {\em osculating hyperplane} of $\cX$ at $P$,
respectively.
   \begin{lemma} {\rm (\cite[Proof of Thm. 1.1]{sv})}\label{lemma2.22} Let
$H$ be a hyperplane in $\P^r(\F)$, and $i\in \{0,\ldots,r-1\}$. Then
$H\supseteq L_i(P)$ if and only if $v_P(\pi^*(H))\ge j_{i+1}(P)$.
   \end{lemma}
In the case where $\F$ is the algebraic closure of a finite field
$\fl$ with $\ell$ elements, and both $\cX$ and $\cD$ are defined
over $\fl$, one can also define the so-called {\em $\fl$-Frobenius
divisor} $S=S_{\cD,\ell}$ associated do $\cD$, see \cite[p. 9]{sv},
whose degree is given by
  \begin{equation}\label{eq2.22}
\deg(S)=\sum_{i=0}^{r-1}\nu_i(2g-2)+(\ell+r)d\, ,
  \end{equation}
where $\nu=0<\ldots<\nu_{r-1}$ is a suitable subsequence of the
$\cD$-order sequence \cite[Prop. 2.1]{sv}.
    \begin{lemma} {\rm (\cite[Prop. 2.4(a), Cor.
2.6]{sv})}\label{lemma2.23} $$v_P(S)\ge
\sum_{i=1}^{r}(j_i(P)-\nu_{i-1})$$ provided that $P\in\cX(\fl)$. In
particular $\cX(\fl)\subseteq\supp(S).$
    \end{lemma}

    \subsection{$\fq$-maximal curves}\label{s2.3} Throughout this
sub-section, $\cX$ denotes an $\fq$-maximal curve of genus $g$. Whenever
concepts and results apply from previous sub-sections, the field $\F$ will
be the algebraic closure $\bar\fq$ of $\fq$. A deep result depending on
the zeta function is the so-called Fundamental Equivalence on divisors
\cite[Cor.1.2]{fgt}:
     \begin{equation}\label{eq2.31}
qP+\frx(P)\sim (q+1)Q\, ,\qquad P\in\cX\, ,\quad Q\in \cX(\fq)\, ,
     \end{equation}
where $\frx$ denotes the Frobenius morphism on $\cX$ relative to
$\fq$. As a consequence, $\cX$ is equipped with the base-point-free
linear series
  $$
\cD_\cX:=|(q+1)P_0|\, , \qquad P_0\in \cX(\fq),
  $$
which is independent of the choice of the point $P_0$ in
$\cX(\fq)$, and has projective dimension $\dim(\cD_\cX)$ at least 2.
Note that (\ref{eq2.31}) is equivalent to
   \begin{equation}\label{eq2.32}
\pi^*(L_{r-1}(P))=qP+\frx(P)\, ,
   \end{equation}
$\pi$ being a morphism associated to $\cD_\cX$. Set $N:=\dim(\cD_\cX)$.
The following result shows that $\cX$ has a non-singular model over $\fq$
given by a curve in $\P^N(\bar\fq)$ of degree $q+1$.
   \begin{lemma} {\rm (Natural embedding theorem, \cite[Thm. 2.5]{kt1},
\cite[Prop.  1.9]{fgt})}\label{lemma2.31} The linear series $\cD_\cX$ is
very ample; i.e. any morphism associated to $\cD_\cX$ is a close
embedding. Equivalently, $q$ is a Weierstrass non-gap at any point of
$\cX.$
   \end{lemma}
The natural embedding theorem together with Castelnuovo's genus
bound (Lemma \ref{lemma2.11}) and its corollary stated in Remark
\ref{rem2.11} provide a very useful upper bound on the genus $g$ of
$\fq$-maximal curves, namely
   \begin{equation}\label{eq2.33}
g\le
\begin{cases}
(q-(N-1)/2)^2/2(N-1) & \text{for odd $N$,}\\
(q-(N-1)/2)^2-1/4)/2(N-1) & \text{for even $N$.}
\end{cases}
   \end{equation}
   \begin{corollary}\label{cor2.31}
\begin{enumerate}
\item[\rm(1)] {\em (\cite{ihara})} $g\le q(q-1)/2;$
\item[\rm(2)] {\em (\cite[Prop. 3]{sti-x})} If $\dim(\cD_\cX)\ge 3$, then
$g\le (q-1)^2/4;$
\end{enumerate}
   \end{corollary}
We point out that Lemma \ref{lemma2.31} together with Corollary
\ref{cor2.31} yields the following lemma that strengthens the
R\"uck-Stichtenoth's characterization of the Hermitian curve \cite{r-sti}
   \begin{lemma} {\rm (\cite[Thm. 2.4]{ft2})}\label{lemma2.32} For a
$\fq$-maximal curve $\cX$ of genus $g$, the following statements are
equivalent:
  \begin{enumerate}
\item[\rm(1)] $g> (q-1)^2/4;$
\item[\rm(2)] $\dim(\cD_\cX)=2;$
\item[\rm(3)] $\cX$ is $\fq$-isomorphic to the Hermitian curve of equation
$Y^qZ+YZ^q=X^{q+1};$
\item[\rm(4)] $g=q(q-1)/2.$
   \end{enumerate}
   \end{lemma}
As a consequence, we have the following result
  \begin{corollary} {\rm (\cite{ft1})}\label{cor2.32} The genus $g$ of a
$\fq$-maximal curve satisfies either
$$
g\le \lfloor (q-1)^2/4\rfloor\qquad\text{or}\qquad g=q(q-1)/2\, .
$$
   \end{corollary}
   \begin{remark}\label{rem2.31} Castelnuovo's number
$c_0(q+1,N)$ in (\ref{castelnuovo}) is attained by an
$\fq$-maximal curve in the cases $q\equiv N-2, 0\pmod{(N-1)}$.
The existence of such a curve $\cX$ is strongly related to the existence
of a point $P_1\in\cX(\fq)$ such that $m$ is a Weierstrass non-gap at $P_1$
satisfying $m(N-1)\le q+1\, (*)$. Since $m_N(P_1)=q+1$ and
$m_{N-1}(P_1)=q$ by Lemma \ref{lemma2.33}(2), if $m$ is a Weierstrass
non-gap at $P_1$, then $m$ must satisfy $m(N-1)\ge q$. Hence, property
$(*)$ occurs when either $m(N-1)=q+1\, (*_1)$ or $m(N-1)=q\, (*_2)$. The
smallest possibilities for $N$ are investigated in the sequel, namely
$N=3$ in Sect. \ref{s3} while $N\in\{4,5\}$ in Sect. \ref{s4}.

In case $(*_1)$, $g=c_0(q+1,N)=(q-1)((q+1)/(N-1)-1)/2$ by \cite[Remark
2.6(1)]{kt1}. There exists just one $\fq$-maximal curve (up to
$\fq$-isomorphism) satisfying $(*_1)$, namely the non-singular $\fq$-model
of the plane curve of equation $y^q+y=x^{(q+1)/(N-1)}$ \cite[Thm.
2.3]{fgt}.

In case $(*_2)$, $g=c_0(q+1,N)=q(q-(N-1))/2(N-1)$ by \cite[Remark
2.6(1)]{kt1}. van der Geer and van der Vlugt, see \cite[Thm.
3.1]{geer-vl2} and \cite[Remark 5.2]{geer-vl3}, by means of fibre product
of certain Artin-Schreier $p$-extensions of the projective line showed
that such curves do exist. Garcia and Stichtenoth, see \cite[Sect. V, Ex.
E]{g-sti1}, noticed that such curves admit a plane model of type
   \begin{equation}\label{eq2.331}
F(y)=f(x)\, ,
   \end{equation}
where $F\in \fq[Y]$ is a $p$-linear polynomial of degree $q/(N-1)$ whose
linear coefficient is different from zero, and where $f\in\fq[X]$ is a
polynomial of degree $q+1$. Here $P_1$ is the unique point over
$x=\infty$. For $N-1=p$, see also \cite[Ex. 1.2]{g-q} and \cite[Prop.
3.5]{geer-vl4}. Unlike the previous case, several pairwise non
$\fq$-isomorphic $\fq$-maximal curves satisfying $(*_2)$ are known to
exist; see \cite[Sect. 5]{abdon}. It has been conjectured \cite[p.
46]{fgt} that a plane $\fq$-model for a $\fq$-maximal curve satisfying
$(*_2)$ has equation of type (\ref{eq2.331}) with $f(x)=x^{q+1}$.
Conversely, the following question arises: Determine the polynomials $F$
and $f$ such that such that the plane curve of equation (\ref{eq2.331})
has an $\fq$-maximal non-singular model. Examples of such curves arise for
instance in \cite{geer-vl-1}, \cite{geer-vl0}, and \cite{geer-vl4}.
Examples of $\fq$-maximal curves defined by (\ref{eq2.331}), where either
$F$ or $f$ are $\fq$-rational functions, can be found in \cite{geer-vl4}
and \cite{g-q}.
   \end{remark}
Finally, some results on Weierstrass Point Theory and Frobenius
orders with respect to the linear series $\cD_\cX$. With the same
notation as in Sect. \ref{s2.2}, Lemma \ref{lemma2.31} together with
(\ref{eq2.31}) forces the first $N$ non-gaps at $P\in\cX$ to have
the following behaviour:
   \begin{equation}\label{eq2.34}
m_1(P)<\ldots<m_{N-1}(P)=q<m_N(P)\, .
   \end{equation}
Furthermore,
   \begin{lemma} {\rm (\cite[Thm. 1.4, Prop.
1.5(ii)(iii)]{fgt})}\label{lemma2.33}
  \begin{enumerate}
\item[\rm(1)] $j_1(P)=1$ for any $P;$ $j_N(P)=q+1$ if $P\in\cX(\fq),$ and
$j_N(P)=q$ otherwise$;$
\item[\rm(2)] $j_{N-i}(P)+m_i(P)=q+1$ for $i=0,\ldots,N,$ provided that
$P\in\cX(\fq);$
\item[\rm(3)] $q-m_i(P)$ is a $(\cD_\cX,P)$-order for $i=0,\ldots,N-1,$
provided that $P\not\in\cX(\fq);$
\item[\rm(4)] $\epsilon_N=\nu_{N-1}=q;$
\item[\rm(5)] $\nu_1=1$ if $N\ge 3.$
  \end{enumerate}
  \end{lemma}
Then, we have one of the main features of the linear series
$\cD_\cX$, namely
$$
\cX(\fq)\subseteq \supp(R_{\cD_\cX})\, .
$$
   \begin{lemma}\label{lemma2.34} Let $\cX$ be a $\fq$-maximal curve of
genus $g$. Set $N:=\dim(\cD_\cX)$.
\begin{enumerate}
  \item[\rm(1)] If $\cX$ is hyperelliptic, then  $q\le 2N-2.$
  \item[\rm(2)] The curve $\cX$ is hyperelliptic provided that
  either
$j_{N-1}(P)=j_N(P)-2$ for $P\in\cX(\fq)$, or $j_{N-1}(P)=j_N(P)-1$
otherwise$.$
  \item[\rm(3)] If there exists $P\in\cX(\fq)$
with $j_{N-1}(P)=j_N(P)-1$, then $q=N-1$.
\end{enumerate}
   \end{lemma}
   \begin{proof} If $\cX$ is hyperelliptic, $m_1(P)=g+1$ at a general
point $P$. Then from (\ref{eq2.34}), $m_{N-1}(P)=g+N-1=q$ and so
$g=q-N+1$. On the other hand $\#\cX(\fq)\le 2(q^2+1)$ and
maximality of $\cX$ yields $2g\le q$. From these computations (1)
follows. Let $P\in\cX(\fq)$ such that $j_{N-1}(P)\in\{q-1,q\}$.
Then from Lemma \ref{lemma2.33}(2) we have $m_1(P)\in\{2,1\}$ and
so either $\cX$ is hyperelliptic or $m_N=N=q+1$. Finally, let
$P\not\in\cX(\fq)$ such that $j_{N-1}(P)=q-1$. Then from
(\ref{eq2.31}), $(q-1)P+D\sim qP+\frx(P)$ with $P\not\in\supp(D)$,
so that $D\sim P+\frx(P)$; i.e. $\cX$ is hyperelliptic.
   \end{proof}
   \begin{lemma}\label{lemma2.35} Let $\cX$ be a $\fq$-maximal curve so
that $j_{N-1}(P)=N-1$ for every point $P\in\cX$, where
$N=\dim(\cD_\cX)$. Then
$$
(N-1)N(g-1)=(q+1)(q-N)\, .
$$
   \end{lemma}
   \begin{proof} The set of $\cD_\cX$-Weierstrass points of $\cX$
   coincides with the set of $\fq$-rational points, and $v_P(R_{\cD_\cX})=1$ for
$P\in\cX(\fq)$; cf. Lemmas \ref{lemma2.33}(1), \ref{lemma2.21}.
Hence the result follows from (\ref{eq2.21}) taking into account
the maximality of $\cX$.
    \end{proof}

    \section{On maximal curves embedded in a quadric surface}\label{s3}

The R\"uck-Stichtenoth theorem together with \cite[Thm. 2.4]{ft2},
stated in the previous section as Lemma \ref{lemma2.32}, gives a
complete classification of $\fq$-maximal curves of genus
$g>(q-1)^2/4$. The objective of this section is to obtain a similar
theorem valid for $(q^2-q+4)/6<g\le (q-1)^2/4$. Notation and
terminology are the same as in Sect. \ref{s2}.
  \begin{theorem}\label{thm3.1} Let $\cX$ be a $\fq$-maximal curve of
genus $g$, and $\pi$ a $\fq$-morphism associated to $\cD_\cX$. Assume
$q\ge 7$. Then the following conditons are equivalent:
   \begin{enumerate}
\item[\rm(1)] $\lfloor(q^2-q+4)/6\rfloor<g\le \lfloor(q-1)^2/4\rfloor;$
\item[\rm(2)] $\dim(\cD_\cX)=3,$ $\pi(\cX)$ lies on a quadric
surface in $\P^3,$ and $g\neq (q^2-2q+3)/6$ whenever $q\equiv 3,5\pmod{6};$
\item[\rm(3)] $\dim(\cD_\cX)=3, \dim(2\cD_\cX)=8,$ and $g\neq
(q^2-2q+3)/6$ whenever $q\equiv 3,5\pmod{6};$
\item[\rm(4)] $\dim(\cD_\cX)=3$ and there exists $P\in \cX(\fq)$ such that
$j_2(P)=(q+1)/2$ if $q$ is odd, or $j_2(P)=(q+2)/2$ otherwise$;$
\item[\rm(5)] $\cX$ is $\fq$-isomorphic to the non-singular $\fq$-model of
either $y^q+y=x^{(q+1)/2}$ if $q$ is odd, or
$y^{q/2}+y^{q/4}+\ldots +y^2+y=x^{q+1}$ otherwise$.$
\item[\rm(6)] $g=(q-1)^2/4$ if $q$ is odd or $g=q(q-2)/4$ otherwise. In
particular the genus $g$ equals Castelnuovo's number $c_0(q+1,3).$
  \end{enumerate}
  \end{theorem}
Under stronger hypotheses, this theorem was partially proved in
\cite[Prop. 2.5]{ft2} for $q$ odd, and in \cite{at} for $q$ even.
  \begin{remark}\label{rem3.1} For $q=2,3,4,5$ the spectrum of the genera
of $\fq$-maximal curves is  $\{0,1\}, \{0,1,3\}, \{0,1,2,6\},
\{0,1,2,3,4,10\}$, respectively; see \cite[Remark 6.1]{g-sti-x}.
  \end{remark}
{}From Theorem \ref{thm3.1} and Remark \ref{rem3.1}, Corollary
\ref{cor2.32}
can be strengthen as follows:
  \begin{corollary}\label{cor3.1} The genus $g$ of a $\fq$-maximal curve
satisfies either
$$
g\le \lfloor(q^2-q+4)/6\rfloor\qquad\text{or}\qquad
g=\lfloor\mbox{$\frac{(q-1)^2}{4}$}\rfloor\qquad\text{or}\qquad
g=(q-1)q/2\, .
$$
   \end{corollary}
   \begin{remark}\label{rem3.2} $\fq$-maximal curves of genus $\lfloor
(q^2-q+4)/6\rfloor$ do exist as the following examples show, see
\cite{g-sti-x}, \cite[Thm. 2.1]{ckt2}:
\begin{enumerate}
\item[\rm(i)] If $q\equiv 2\pmod{3}$, the non-singular $\fq$-model of
$x^{(q+1)/3)}+x^{2(q+1)/3}+y^{q+1}=0$ is $\fq$-maximal and has genus
$(q^2-q+4)/6$.

\item[\rm(ii)] If $q\equiv 1\pmod{3}$, the non-singular $\fq$-model of
$y^q-yx^{2(q-1)/3}+x^{(q-1)/3}=0$ is $\fq$-maximal and has genus
$(q^2-q)/6$.

\item[\rm(iii)] If $q=p^t\equiv 0\pmod{3}$, the non-singular $\fq$-model
of $y^q+y+(\sum_{i=1}^{t}x^{q/p^i})^2=0$ is $\fq$-maximal and has genus
$(q^2-q)/6$.
   \end{enumerate}
It may be that no further infinite family exists. Also, each of the
above curves is $\fq$-covered by the Hermitian curve via a suitable
morphism of degree $3,$ and it would be of interest to prove or
disprove uniqueness of some (perhaps all) of these examples.
   \end{remark}
   \begin{remark}\label{rem3.3} In searching quantitative results for the
number of $\fl$-rational points of a curve of genus $g$, the maximum
number $N_\ell(g)$ of $\fl$-rational points on such curves play an
important role; see e.g. \cite{geer-vl2}. Corollary \ref{cor3.1} excludes
certain values for $N_{q^2}(g)$ whenever $(q^2-q+4)/6<g<(q-1)^2/4$ or
$(q-1)^2/4<g<q(q-1)/2$. More precisely, for such values of $g$, we have
$N_{q^2}(g)<q^2+1+2qg$. A similar resuly follows from Theorem
\ref{thm4.11}(a). Hence from deeper results due to J.P. Serre and K.
Lauter one can deduce $N_{q^2}(g)\le q^2+1+2qg-m$, where $m\in\{1,2,3\}$,
cf. \cite{geer-vl}. One can also obtain improvements on some entries in
the tables of loc. cit. For instance, we have $N_{64}(11)\leq 238$,
$N_{81}(13)\leq 314$, $N_{81}(15)\le 350$, while the upper bounds in the
tables are respectively 241, 316, 352. It should be noted that the above
considerations will extend to a more general case, once the conjecture
stated in the introduction has been proved. 
   \end{remark}
In proving Theorem \ref{thm3.1}, we will need some technical
results concerning $\fq$-maximal $\cX$ with $\dim(\cD_\cX)=3$.
  \begin{lemma}\label{lemma3.1} Let $\cX$ be a $\fq$-maximal curve with
$\dim(\cD_\cX)=3$, and $\pi$ a $\fq$-morphism associated to $\cD_\cX$.
Assume $q\ge 4.$
\begin{enumerate}
\item[\rm(1)] $\dim(2\cD_\cX)\ge 8$.
\item[\rm(2)] If $\dim(2\cD_\cX)=8$, then $\pi(\cX)$ lies on
a quadric surface in $\P^3(\bar\fq)$.
\item[\rm(3)] The quadric surface $\cQ$ in part (2) is uniquely determinated by
the property $\pi(\cX)\subseteq \cQ$, and it is defined over $\fq.$
\end{enumerate}
  \end{lemma}
  \begin{proof} (1) Let $P\in \cX(\fq)$ and set $m_i:=m_i(P)$. From Lemma
\ref{lemma2.33}(2), $m_2=q$ and $m_3=q+1$. Then, as $2m_1\ge m_2=
q$ and $q\ge 4$, it is easy to see that there are at
least 8 positive Weierstrass non-gaps in $[m_1,2m_3]$
and so $\dim(2\cD_\cX)\ge 8$.

(2) See \cite[p. 352]{hartshorne}.

(3) If $\pi(\cX)$ lies on $\cQ$, then $\pi(\cX)$ also lies on
${\mathbf Fr}(\cQ)$, where ${\mathbf Fr}$ is the Frobenius
collination on $\P^3(\bar\fq)$ relative to $\fq$. Clearly
$\cQ={\mathbf Fr}(\cQ)$ if and only if $\cQ$ is defined over $\fq$.
It this were not the case in our situation, then $\cX$ would be
contained in the intersection of two distinct quadrics,
contradicting the hypothesis $q+1=\deg(\pi(\cX))\le 4$ by the
B\'ezout theorem.
  \end{proof}
  \begin{lemma}\label{lemma3.2} Let $\cX$ be a $\fq$-maximal curve with
$\dim(\cD_\cX)=3$, $\pi$ a morphism associated to $\cD_\cX$, and
$P\in\cX$. Suppose that
$\pi(\cX)$ lies on a quadric surface $\cQ$ in $\P^3(\bar\fq)$, and that
$q\ge 5$. Then
  \begin{enumerate}
\item[\rm(1)] $j_2(P)\in\{2,j_3(P)/2,(j_3(P)+1)/2\};$
\item[\rm(2)] $j_2(P)>2$ if and only if the tangent line $L_1(P)$
of $\cX$ at $P$ lies on $\cQ;$
\item[\rm(3)] either $q$ is even, $j_2(P)=q/2$ and $P\not\in\cX(\fq)$ or
$q$ is odd, $j_2(P)=(q+1)/2$ and $P\in\cX(\fq)$ provided that
$j_2(P)>2$ and that $\cQ$ is non-singular at $\pi(P).$
  \end{enumerate}
\end{lemma}
   \begin{proof} Set $j_i:=j_i(P)$, $i=0,\ldots,3$. Let $x_0=1,
x_1,x_2,x_3$ be $\fq$-rational functions on $\cX$, such that
$v_P(x_i)=j_i$. Up to a projective collineation in $\P^3(\bar\fq)$,
we can assume $\pi=(x_0:x_1:x_2:x_3)$. Let $(X_0,\ldots,X_3$) be
coordinates in $\P^3(\bar\fq)$ such that each $x_i$ is the
pull-back via $\pi$ of $X_i/X_0$ restricted to $\pi(\cX)$. Then
$\pi(P)=(1:0:0:0)$ and $L_1(P)$ is given by $X_2=X_3=0$; see
\cite[proof of Thm. 1.1]{sv}. Let the quadric $\cQ$ have
homogeneous equation
  \begin{align*}
F(X_0,X_1,X_2,X_3)=\, \, &
a_{00}X_0^2+a_{01}X_0X_1+a_{02}X_0X_2+a_{03}X_0X_3+
a_{11}X_1^2+a_{12}X_1X_2+\\
  \, \, &  a_{13}X_1X_3+a_{22}X_2^2+a_{23}X_2X_3+a_{33}X_3^2\, .
  \end{align*}
Then $a_{00}=0$ because of $F(\pi(P))=0$. Furthermore, $x_1,x_2$
and $x_3$ are related in the function field over $\bar\fq$ of $\cX$ by
$F(1,x_1,x_2,x_3)=0$. In
addition, the valuation at $P$ of the functions $x_1, x_2, x_3,
x_1^2, x_1x_2, x_1x_3, x_2^2, x_2x_3, x_3^2$ are respectively
   \begin{equation}\label{eq3.1}
1,j_2,j_3, 2, j_2+1, j_3+1, 2j_2, j_3+j_2, 2j_3\, .
   \end{equation}
Hence, $a_{01}=0$.

(1) $j_2+1<j_3$ by Lemma \ref{lemma2.34} and the hypothesis $q\ge
5$. So from the inequalities
$$
2\le j_2<j_2+1<j_3<j_3+1<j_3+j_2<2j_3
$$
and (\ref{eq3.1}) we obtain part (1).

(2) We have from (\ref{eq3.1}) that $j_2>2$ if and only if
$a_{11}=0$. Now, as $F(X_0,X_1,0,0)=a_{11}X_1^2$, the last
condition is equivalent to $L_1(P)\subseteq\cQ$ and the result
follows.

(3) If $j_2>2$, from the proof of part (1) we get
$a_{11}=a_{02}=a_{12}=0$. An easy computation shows then that $\cQ$
is non-singular at $\pi(P)$ if and only if $a_{03}\neq 0$.
Therefore $2j_2=j_3$, and the result follows from Lemma
\ref{lemma2.33}(1).
   \end{proof}
    \begin{proposition}\label{prop3.1} Let $\cX$ be a $\fq$-maximal curve
and $\pi$ a $\fq$-morphism associated to $\cD_\cX$. Suppose that
$q$ is even, $q>4$, and that $\pi(\cX)$ lies on a quadric $\cQ$ in
$\P^3(\bar\fq)$. Then
  \begin{enumerate}
\item[\rm(1)] $\cQ$ is a cone$;$
\item[\rm(2)] the vertex $V$ of $\cQ$ belongs to $\pi(\cX)$; if
$V=\pi(\tilde V)$, then $\tilde V\in\cX(\fq)$ and $j_2(\tilde V)=(q+2)/2.$
   \end{enumerate}
\end{proposition}
\begin{proof} General properties of quadrics of a $3$-dimensional projective space
over a finite field can be found in \cite{hirschfeld}. Here we will use
the following properties:
Let $P\in\cQ$ be a non-singular point of
$\cQ$ and denote by $T_P\cQ$ the tangent plane of $\cQ$ at $P$.
  \begin{itemize}
\item If $P\in\pi(\cX)$, then $T_P\cQ\supseteq L_1(P)$;
\item Let $\ell$ and $\ell_1$ be lines such that $P\in\ell\subseteq\cQ$, and
$\ell_1\subseteq T_P\cQ$. If $\ell\neq\ell_1$, then $T_P\cQ$ is generated by
$\ell$ and $\ell_1$;
\item There exist lines $\ell$ and $\ell_1$ such that
$P\in\ell\cap\ell_1$, and $\cQ\cap T_P\cQ=\ell\cup\ell_1$;
  \end{itemize}
If $\cQ$ is non-singular, then
  \begin{itemize}
\item No two tangent hyperplanes of $\cQ$ at different points coincide.
  \end{itemize}

To simplify our notation we shall identify $\cX$ and $\pi(\cX)$,
according to Lemma \ref{lemma2.31}.

(1) Since $\cX$ is non-degenerate, $\cQ$ is irreducible. Then $\cQ$ is a
cone if and only if $\cQ$ is singular, as this case can only occur when
$\cQ$ has just one singular point.

Suppose that $\cQ$ is non-singular. Then from Lemma
\ref{lemma3.2}(3), $j_2(Q)=2$ for each $Q\in\cX(\fq)$. Note that
there exists $P\in\cX\setminus\cX(\fq)$ such that $j_2(P)>2$; in
fact, otherwise Lemma \ref{lemma2.35} would yield
$6(g-1)=(q+1)(q-3)$; but then $q$ would be odd, a contradiction.
Hence $j_2(P)=q/2$ by Lemma \ref{lemma3.2}(3). Let $Q_1\in
\cX(\fq)$. We have $Q_1\not\in L_1(P)$, as $\cX\cap
L_1(P)\subseteq\cX\cap L_2(P)=\{P,\frx(P)\}$ (cf. (\ref{eq2.32})),
and hence the plane $H=H_{Q_1}$ generated by $L_1(P)$ and $Q_1$ is
well defined. Then $H\neq L_2(P)$, and the intersection divisor of
$\cX$ and $H$ becomes
   \begin{equation}\label{eq3.2}
\cX\cdot H=\frac{q}{2}P+D\, ,
   \end{equation}
where $D=D_{Q_1}$ is a divisor on $\cX$ of degree $(q+2)/2$ with
$Q_1\in\supp(D)$, and $P\not\in\supp(D)$. In addition, Lemma
\ref{lemma3.2}(2) assures the existence of a line $\ell=\ell_{Q_1}$
such that
   \begin{equation}\label{eq3.3}
\cQ\cap H=L_1(P)\cup\ell\, .
   \end{equation}
Actually, the line $\ell$ is defined over $\fq$. In fact, $\cQ$ is
defined over $\fq$ by Lemma \ref{lemma3.1}(3), and
$Q_1\in\cX(\fq)\setminus L_1(P)$ implies that $Q_1\in \ell$.
  \begin{claim}\label{claim3.1} $\cX\cap\ell\subseteq \cX(\fq).$
  \end{claim}
{\em Proof of Claim \ref{claim3.1}.} If there exists
$Q\in\cX\cap\ell\setminus \cX(\fq)$, then $\frx(Q)\in\ell$ as
$\ell$ is defined over $\fq$. Thus $\ell\subseteq L_2(Q)$, and
hence $\ell\cap\cX\subseteq\{Q,\frx(Q)\}$. It follows
$Q_1\not\in\ell$, but this is a contradiction.
  \begin{claim}\label{claim3.2} If
$Q\in\supp(D)\setminus\{\frx(P)\}$, then $Q\in\cX(\fq)$ and
$v_Q(D)=1.$
  \end{claim}
{\em Proof of Claim \ref{claim3.2}.} Since
$\supp(D)\setminus\{\frx(P)\}\subseteq \ell\cap\cX$, we have
$Q\in\cX(\fq)$ by Claim \ref{claim3.1}. Now if $v_Q(D)\ge 2$, then
$H\supseteq L_1(Q)$ by $j_2(Q)=2$ and Lemma \ref{lemma2.22}. Also,
$\ell\neq L_1(Q)$ because $L_1(Q)\not\subseteq\cQ$ by Lemma
\ref{lemma3.2}(2). Therefore the plane $H$ is generated by the
lines $\ell$ and $L_1(Q)$, and hence $H=T_{Q_1}\cQ$. Let $\ell_1$
be the line defined by $Q_1\in\ell_1$, and $\cQ\cap
T_{Q_1}\cQ=\ell\cap\ell_1$. From (\ref{eq3.3}), we infer that
$L_1(P)=\ell_1$ and so $Q_1\in L_1(P)$, but this is a contradiction.
   \begin{claim}\label{claim3.3} $\frx(P)\not\in\supp(D).$
  \end{claim}
{\em Proof of Claim \ref{claim3.3}.} Suppose on the contrary that
$\frx(P)\in\supp(D)$. Equivalently, $\frx(P)\in L_1(P)$ by Claim
\ref{claim3.1}. Then $v_{\frx(P)}(D)=1$. In fact, using a similar
argument to that in the proof the previous claim, one can show that
$v_{\frx(P)}(D)\neq 1$ together with $L_1(P)\neq L_1(\frx(P))$
implies $H=T_{\frx(P)}\cQ$ in contradiction with (\ref{eq2.32}).
Hence, for each $Q\in \cX(\fq)$, the divisor $D$ in (\ref{eq3.2})
may also be written as $D=D_Q=\frx(P)+D_Q'$ in such a way that
(\ref{eq3.3}) holds true, $\supp(D_Q')
\subseteq\cX(\fq)$, and $\deg(D_Q')=q/2$. Notice that $H_Q$ is
generated by $L_1(P)$ and $Q'$ $(*)$ where $Q'$ is any point of
$\supp(D_Q')$. Now let $Q_1, Q_2\in\cX(\fq)$ such that $Q_2\not\in
\supp(D_{Q_1}')$. Then $\supp(D_{Q_1}')\cap\supp(D_{Q_2}')=\emptyset$,
otherwise $H_{Q_1}=H_{Q_2}$ by $(*)$. This yields that $q/2$ must
divide the number of $\fq$-rational points of $\cX$, which is a
contradiction because $\#\cX(\fq)=q^2+1+2gq$ is an odd number.

So far we have shown that each $Q_1\in\cX(\fq)$ gives rise to a
plane $H_{Q_1}$, to a line $\ell=\ell_{Q_1}$, and to a divisor
$D=D_{Q_1}$ such that (\ref{eq3.2}) and (\ref{eq3.3}) hold with
$D=Q_1+Q_2+\ldots+Q_{(q+2)/2}$ being the sum of $(q+2)/2$
$\fq$-rational points. Notice that $\supp(D)=\cX\cap\ell$. Let
$\ell_1$ be chosen in such a way that $Q_1\in
\ell_1$ and that
  \begin{equation}\label{eq3.4}
\cQ\cap T_{Q_1}\cQ=\ell\cup\ell_1\, .
  \end{equation}
Clearly, $\ell_1$ is $\fq$-rational, and thus
$\cX\cap\ell_1\subseteq
\cX(\fq)$ as in the proof of Claim
\ref{claim3.1}. Therefore
\begin{equation}\label{eq3.5}
\cX\cdot T_{Q_1}\cQ=2Q_1+Q_2+\ldots+Q_{(q+2)/2}+D'\, ,
\end{equation}
where $D'$ is a divisor on $\cX$ of degree $(q-2)/2$ such that
$Q_1\not\in \supp(D')\subseteq \cX(\fq)$.
   \begin{claim}\label{claim3.4} $\supp(D)\cap\supp(D')=\emptyset$, and
$v_S(D')=1$ for each $S\in \supp(D').$
  \end{claim}
{\em Proof of Claim \ref{claim3.4}.} Let $S\in\supp(D')$. Suppose
on the contrary that $S=Q_i$ for some $i$. Then $T_{Q_1}\cQ$
contains $L_1(Q_i)$ which is different from $\ell$ as $j_2(Q_i)=2$.
Hence $T_{Q_1}\cQ$ is generated by $L_1(Q_i)$ and $\ell$. These
lines also generate $T_{Q_i}\cQ$ and so $i=1$ contradicting
$Q_1\not\in\supp(D')$.

Finally suppose on the contrary that $v_S(D_2)\ge 2$. Replacing
$\ell$ by $\ell_1$, the above argument shows that
$T_S\cQ=T_{Q_1}\cQ$, whence $S=Q_1$ follows, again a contradiction.

Therefore, to each $Q_1$ we have associated two lines $\ell$ and
$\ell_1$ such that both (\ref{eq3.4}) and (\ref{eq3.5}) hold where
$D'$ is a divisor of degree $(q-2)/2$, $\supp(D')\subseteq\cX(\fq)$,
and $\supp(D)\cap\supp(D')=\{Q_1\}$. As it is well-known, $\cQ$ has
just two families of lines contained in $\cQ$ and any two lines of
the same family are disjoint. This implies again that $\#\cX(\fq)$
must be a multiple of $q/2$, contradicting the $\fq$-maximality of
$\cX$.

(2) As $q$ is even, from Lemma \ref{lemma2.35} there exists
$P\in\cX$ such that $j_2(P)>2$. Suppose that $P\not\in\cX(\fq)$.
{}From $j_2(P)P+D\sim (q+1)P_0$, we find that
$j_2(P)\frx(P)+\frx(D)\sim (q+1)P_0$ and so
$j_2(\frx(P))=j_2(P)>2$. Therefore $L_1(P)\cup L_1(\frx(P)\subseteq
\cQ$ by Lemma \ref{lemma3.2}(2), and hence $V\in L_1(P)\cap
L_1(\frx(P))$. Now, since $V$ is $\fq$-rational by Lemma
\ref{lemma3.1}(3), we have $\frx(P)\neq V$, and hence
$L_1(\frx(P))$ is generated by $\frx(P)$ and $V$; in particular
$L_1(\frx(P))\subseteq L_2(P)$ and thus $1=v_{\frx(P)}(\cX\cdot
L_2(P))\ge j_2(\frx(P))$ by Lemma \ref{lemma2.22}, a contradiction.

Therefore $P$ must be $\fq$-rational and hence $\cQ$ must have a
singularity at $P$ by Lemma \ref{lemma3.2}(3). Then $P=V$ and
$j_2(P)=(q+2)/2$ by Lemma \ref{lemma3.2}(1) and the assumption of $q$
being even.
  \end{proof}

{\em Proof of Theorem \ref{thm3.1}.} $\text{(1)$\Rightarrow$(2)}:$
{}From the hypothesis on $g$, $\dim(\cD_\cX)=3$ follows by
(\ref{eq2.33}) and Lemma \ref{lemma2.32}. Since $c_1(q+1,3)$ in
Lemma \ref{lemma2.13} is equal to $\lfloor (q^2-q+4)/6\rfloor$,
that lemma together with Lemma \ref{lemma2.31} shows that
$\pi(\cX)$ lies on a quadric provided that
$q\not\in\{7,8,9,11,13,17,19,23\}$.

Assume $q=8$. Then $g>(q^2-q+4)/6=10$. By virtue of Lemma
\ref{lemma3.1}(1)(2), it is enough to show that $\dim(2\cD_\cX)\le
8$. Suppose on the contrary that $\dim(2\cD_\cX)\ge 9$. Then from
Lemma \ref{lemma2.11} and Remark \ref{rem2.11}, $g\le
(q-1)(q-2)/4=10.5$ follows, a contradiction.

Now, let $q$ be odd, $q\ge 7$. Our goal is to show that the second
positive $\cD_\cX$-order $\epsilon_2$ (see sections \ref{s2.2},
\ref{s2.3}) is equal to two. In fact, if this is the case, then the
Generic Order of Contact Theorem \cite[Thm. 3.5]{hefez-kleiman}
yields that the curve $\cX$ (that is $\pi(\cX)$ by previous
identification) is reflexive. Reflexivity forces the monodromy
group of $\cX$ to be isomorphic to the symmetric group $S_{q+1}$,
see (\cite[p. 264]{ballico-hefez}, \cite[Cor. 2.2]{rathmann}).
Hence the points of a general hyperplane section of $\cX$ lie in
uniform position \cite[Cor. 1.8]{rathmann}. Then Lemma
\ref{lemma2.13} holds true; see Remark \ref{rem2.12}.

Suppose on the contrary that $\epsilon_2>2$. Let $S$ be the
$\fq$-Frobenius divisor associated to $\cD_\cX$. From Lemmas
\ref{lemma2.23}, \ref{lemma2.21}(1), \ref{lemma2.33}(4)(5),
$v_P(S)\ge
\epsilon_2+1\ge 4$ for any $P\in\cX(\fq)$. Then by (\ref{eq2.22})
and the $\fq$-maximality of $\cX$, $(3q-1)(2g-2)\le
(q+1)(q^2-4q-1)$. On the other hand, $2g-2>(q+1)(q-2)/3$ by
hypothesis, and thus $5q+5<0$, a contradiction.

$\text{(3)$\Rightarrow$(2)}:$ This follows from Lemma \ref{lemma3.1}(2).

$\text{(2)$\Rightarrow$(4)}:$ Let $q$ be odd. There exists
$P\in\cX$ such that $j_2(P)>2$, otherwise $g$ would be equal to
$(q^2-2q+3)/6$ by Lemma \ref{lemma2.35}. If such a point $P\in\cX$
should not be in $\cX(\fq)$, then by Lemma \ref{lemma3.2}(3) both
$P$ and $\frx(P)$ would be singular points of the quadric, a
contradiction. Therefore $P\in\cX(\fq)$ and hence $j_2(P)=(q+1)/2$
by Lemma \ref{lemma3.2}(1). If $q$ is even, the result follows from
Proposition \ref{prop3.1}(2).

$\text{(4)$\Rightarrow$(5)}:$ From Lemma \ref{lemma2.33}(2) and the
hypothesis, $m_1(P)=(q+1)/2$ for $q$ is odd, and $m_1(P)=q/2$ for
$q$ even. In the odd case, $(\dim(\cD_\cX)-1)m_1(P)=q+1$, and (5)
follows from \cite[Thm. 2.3]{fgt}. In the even case,
$(\dim(\cD_\cX)-1)m_1(P)=q$, and hence $g=q(q-2)/4$ by \cite[Remark
2.6(1)]{kt1}. Then (5) follows from the main result in \cite{at}.

The implications $\text{(5)$\Rightarrow$(6)}$,
$\text{(6)$\Rightarrow$(1)}$, and $\text{(5)$\Rightarrow$(3)}$ are
trivial.

\section{On $\fq$-maximal curves whose genus equals Castelnuovo's
number}\label{s4}

In this section we investigate certain $\fq$-maximal curves whose genus
equals Castelnuovo's number $c_0(q+1,N)$ for $N\in\{4,5\}$.

\subsection{The case $q\equiv1,2\pmod{3}$}\label{s4.1} The main result is
Theorem \ref{thm4.11} which provides a complete description of
$\fq$-maximal curves of genus $g=(q-1)(q-2)/6$, $q\equiv
1,2\pmod{3}$, $q\ge 11$: Such $\fq$-maximal curves can only exist for
$q\equiv 2\pmod{3}$, and they are $\fq$-isomorphic to the non-singular
$\fq$-model of the plane curve of equation
   \begin{equation}\label{eq4.11}
y^q+y=x^{(q+1)/3}\, .
   \end{equation}
To do this let $\cX$ denote an $\fq$-maximal curve of genus
$g=(q-1)(q-2)/6$ with $q\equiv 1,2\pmod{3}$, equipped with the
linear series $\cD_\cX$ as defined before. The first step is to
compute the dimension of $\cD_\cX$.
  \begin{lemma}\label{lemma4.11} $\dim(\cD_\cX)=4.$ In particular,
$g=c_0(q+1,4).$
  \end{lemma}
  \begin{proof} From (\ref{eq2.33}) and Lemma \ref{lemma2.32},
$\dim(\cD_\cX)\in\{3,4\}$. Suppose on the contrary that
$\dim(\cD_\cX)=3$. If $\epsilon_2=2$, (\ref{eq2.21}) becomes $\deg
(R)=(3+q)(2g-2)+4(q+1)$, while $\fq$-maximality of $\cX$ implies
$\deg(R)\geq q^2+1+2gq$ as $v_P(R)\geq 1$ for every $P\in\fq(\cX)$.
But then $g\ge (q^2-2q+3)/6$ contradicting the hypothesis on $g$.
If $\epsilon_2>2$, then $\epsilon_2\ge 5$ by the $p$-adic criteriom
\cite[Cor. 1.9]{sv} and $q\not\equiv 0\pmod{3}$. Replacing the
ramification divisor $R$ by the Frobenius divisor $S$ in the
previous argument yields again a contradiction. In fact,
(\ref{eq2.22}) reads currently $\deg(S)=(1+q)(2g-2)+(q^2+3)(q+1)$,
while $\deg(S) \geq (q^2+1+2gq)(\epsilon_2+1)$ by the
$\fq$-maximality of $\cX$ and the lower bound $v_P(S)\ge
\epsilon_2+1$ for $P\in\fq(\cX)$ which has been shown in the proof
of Theorem \ref{thm3.1}. Taking $\epsilon_2\geq 5$ into account,
this gives $(5q-1)(2g-2)\le (q+1)(q^2-6q-3)$, whence $2q^2-3q+13\le
0$ follows for $g=(q-1)(q-2)/6$; a contradiction.
   \end{proof}
We take advantage of the current hypothesis that the genus of $\cX$
is equal to Castelnuovo's number $c_0(q+1,4)$ by means of
Lemma \ref{lemma2.12}(1). Indeed, this lemma implies that
$\dim(2\cD_\cX)=11$ which allows to compute the possibilities for
$(\cD_\cX,P)$-orders. To show how to do this, set $j_i=j_i(P)$ and
denote by $\Sigma_P$ the set of $(2\cD_\cX,P)$-orders. Then
$\Sigma_P$ contains both the following sets $\Sigma_1$ and
$\Sigma_2$:
  \begin{align}
\Sigma_1:= & \{0,1,2,j_3,j_4,j_4+1, j_4+j_2, j_4+j_3,
2j_4\}\label{eq4.12}\\
\Sigma_2:= & \{j_2,j_2+1,j_3+1,2j_2,j_3+j_2,2j_3\}\, ,\notag
  \end{align}
where $j_4=q+1$ for $P\in\cX(\fq)$, and $j_4=q$ otherwise (cf.
Lemma \ref{lemma2.33}(1)).
  \begin{lemma}\label{lemma4.12} Let $\cX$ be a $\fq$-maximal curve
and $P\in\cX$ a point with $j_2(P)=2$. If $\dim(\cD_\cX)=4$,
$\dim(2\cD_\cX)=11$, and $q\ge 9$, then $j_3(P)=3$.
  \end{lemma}
  \begin{proof} The hypothesis on $q$ together with
Lemma \ref{lemma2.34} implies that
  \begin{equation}\label{eq4.13}
\text{$j_3<j_4-2$ for $P\in\cX(\fq)$\quad and \quad
$j_3<j_4-1$ otherwise}\, .
  \end{equation}
Suppose $j_3>3$. If $P\in \cX(\fq)$, from (\ref{eq4.12}) and
(\ref{eq4.13})
$$
\Sigma_P=\Sigma_1\cup \{3,j_3+1,j_3+2\}\, ,
$$
and $2j_2,2j_3\in\Sigma_P$. Thus $j_3=2j_2=4$ so that
$2j_3=8=j_4=q+1$; i.e. $q=7$. If $P\not\in\cX(\fq)$ and $j_3>4$, from
(\ref{eq4.12}) and (\ref{eq4.13}) we have
$$
\Sigma_P=\Sigma_1\cup\{3,4, j_3+1\}\, ,
$$
and $(j_3+2,2j_3)\in\{(q,q+1), (q,q+2), (q+1,q+2)$. Then $j_3\le
4$, a contradiction. Finally, if $P\not\in\cX(\fq)$ and $j_3=4$,
then (\ref{eq4.12}) together with (\ref{eq4.13}) gives
$$
\Sigma_P=\Sigma_1\cup\{3,5,6,8\}\, .$$
Hence $j_4=q=8$, and this completes the proof.
   \end{proof}
   The previous lemma together with Lemma
\ref{lemma2.35} gives the following result.
   \begin{corollary}\label{cor4.11} Let $\cX$ be a $\fq$-maximal curve
such that $\dim(\cD_\cX)=4$ and $\dim(2\cD_\cX)=11$. Assume $q\ge 9$. If
$j_2(P)=2$ for any $P\in\cX,$ then $q\equiv 1,2\pmod{3}$
and $g=(q^2-3q+8)/12.$
   \end{corollary}
Now, we investigate the case $j_2(P)>2$ for some $P\in\cX$.

   \begin{lemma}\label{lemma4.13} Let $\cX$ be a $\fq$-maximal curve and
$P\in\cX$ a point with $j_2(P)>2$. Suppose that $\dim(\cD_\cX)=4,$
$\dim(2\cD_\cX)=11,$ and that $q\ge 7$.
   \begin{enumerate}
\item[\rm(1)] If $P\in \cX(\fq)$ and $g>(q-2)q/8$ for $q$ even, then
either $q\equiv 2\pmod{3}, j_2(P)=(q+1)/3, j_3(P)=(2q+2)/3$; or $q\equiv
0\pmod{3}, j_2(P)=(q+3)/3, j_3(P)=(2q+3)/3;$
\item[\rm(2)] If $P\not\in\cX(\fq)$, then either $q\equiv 1\pmod{3},
j_2(P)=(q+2)/3, j_3(P)=(2q+1)/3;$ or $q\equiv 0\pmod{3},
j_2(P)=q/3, j_3(P)=2q/3;$ or $q$ is odd, $j_2(P)=(q-1)/2,$
$j_3(P)=(q+1)/2;$ or $q$ is even, $j_2(P)=q/2, j_3(P)=(q+2)/2.$
   \end{enumerate}
   \end{lemma}
   \begin{proof} Suppose first that $j_3>j_2+1$. According to (\ref{eq4.12}) and
(\ref{eq4.13}) we have only three possibilities, namely
$$
\Sigma_P=\Sigma_1\cup\{j_2,j_2+1,j_3+1\}\, ,
$$
and $(j_3+j_2,2j_3)\in\{(j_4,j_4+1), (j_4, j_4+j_2),
(j_4+1,j_4+j_2)\}$. The first one cannot actually occur by $j_3\neq
j_2+1$; from the second one $j_4\equiv 0\pmod{3}$, $j_2=j_4/3$,
$j_3=2j_4/3$ follow, while the third one gives $j_4\equiv
1\pmod{3}$, $j_2=(j_4+2)/3$, and $j_3=(2j_4+1)/3$.

Suppose next that $j_3=j_2+1$. Then $2j_2\not\in\{j_3,j_3+1\}$ by
$j_2>2$. Moreover, $2j_2\neq j_4+1$; otherwise $j_2=(j_4+1)/2,
j_3=(j_4+3)/2$ and from (\ref{eq4.12}) and (\ref{eq4.13}) we would
have
$$
\Sigma_P=\Sigma_1\cup\{j_2,j_3+1,j_4+2,j_4+3\}
$$
which implies $j_4+j_2=j_4+3$; whence $j_4=5$ and so $q\le 5$. If
$2j_2=j_4$, then $P\not\in\cX(\fq)$; otherwise $j_3=(q+3)/2$ and
hence $m_1=(q-1)/2$ by Lemma \ref{lemma2.33}(2), and this would
imply $\dim(\cD_\cX)\ge 5$. Finally, assume that
$2j_2\not\in\{j_3,j_3+1,j_4,j_4+1\}$. Then from (\ref{eq4.12}) and
(\ref{eq4.13})
$$
\Sigma_P=\{j_2, j_3+1, 2j_2 \}\, ,
$$
and $j_3+j_2\in\{j_4,j_4+1\}$. If $j_3+j_2=j_4+1$, then $2j_2=j_4$,
whence $j_3+j_2=j_4$. Then $j_2=(j_4-1)/2$ and $j_3=(j_4+1)/2$. We
claim that $P\not\in\cX(\fq)$. In fact, otherwise $j_2=q/2$,
$j_3=(q+2)/2$ and hence $m_1=q/2$, $m_2=(q+2)/2$ by Lemma
\ref{lemma2.33}(2) which yields $g\le (q-2)q/8$, a contradiction.
    \end{proof}
   \begin{theorem}\label{thm4.11} Assume $q\ge 11$.
  \begin{enumerate}
\item[\rm(1)] If $q\equiv 1\pmod{3},$ there is no $\fq$-maximal curve
of genus $(q-1)(q-2)/6.$
\item[\rm(2)] If $q\equiv 2\pmod{3},$ the following statements are
equivalent for a $\fq$-maximal curve $\cX$ of genus $g$:
 \subitem\rm(a) $g=(q-1)(q-2)/6;$
 \subitem\rm(b) $\exists P\in\cX(\fq), \exists m\in H(P)$ such that
$3m=q+1;$
 \subitem\rm(c) $\cX$ is $\fq$-isomorphic to the non-singular $\fq$-model
of the curve (\ref{eq4.11})$.$
  \end{enumerate}
   \end{theorem}
\begin{proof} (1) Suppose on the contrary that $\cX$ is an $\fq$-maximal curve
of genus $g=(q-1)(q-2)/3$ with $q\equiv 1\pmod{3}$. Since
$q+1=\frac{q-1}{3}\cdot 3+2$, we have $g=c_0(q+1,3)$ by Lemma
\ref{lemma4.11}. Hence, Lemma \ref{lemma2.12} implies that
$\dim(2\cD_\cX)=11$ and that $\frac{q-4}{3}\cD_\cX$ is the
canonical linear series on $\cX$. Then
  \begin{equation}\label{eq4.14}
a_1i_1+\ldots+a_{(q-4)/3}i_{(q-4)/3}+1\not\in H(P)\, ,
  \end{equation}
where the $i_j$'s are $(\cD_\cX,P)$-orders, and the $a_j$'s are
non-negative integers such that $\sum_j a_j\le (q-4)/3$. We choose
then $P\in\cX$ with $j_2(P)>2$ according to Corollary
\ref{cor4.11}. By Lemma \ref{lemma4.13}, $P\not\in\cX(\fq)$. Thus,
we have to analyze three cases. As before, $m_i=m_i(P)$ stands for
the $i$th Weierstrass non-gap at $P$. Recall that $m_3=q$ by
(\ref{eq2.34})).

Case 1: $j_2(P)=(q+2)/3, j_3(P)=(2q+1)/3$. From Lemma
\ref{lemma2.33}(3), $\{q-m_2,q-m_1\}\subseteq
\{1,(q+2)/3,(2q+1)/3\}$. We have that $q-m_1=(q+2)/3$, since otherwise
$m_1=(q-1)/3$ and hence $q\ge m_4$, a contradiction. Thus
$m_1=(2q-2)/3$. However this leads again to a contradiction since,
by (\ref{eq4.14}), $(q-7)/3+(q+2)/3+1=(2q-2)/3$ does not belong to
$H(P)$.

Case 2: $q$ odd, $j_2(P)=(q-1)/2, j_3(P)=(q+1)/2$. From
(\ref{eq4.14}), $2j_2(P)+1=q$ does not belong to $H(P)$, a
contradiction.

Case 3: $q$ even, $j_2(P)=q/2, j_3(P)=(q+2)/2$. Arguing as in Case 1
we have either $m_1=q/2-1$ or $m_1=q/2$. In the former case,
$q-2\in H(P)$ and thus Lemma \ref{lemma2.33}(3) implies $j_2(P)=2$.
Since this is not admitted currently, the latter case can only
occur. Then $m_1=q/2$ and $m_2=q-1$. Now, as $\dim(2\cD_\cX)=11$,
from (\ref{eq2.31}) $m_9=2q$ follows. Since a similar result to
Lemma \ref{lemma2.33}(3) holds, namely $2q-m_i$ is a
$(2\cD_\cX,P)$-order for $i=0,\dots,9$, and the set of
$(2\cD_\cX,P)$-orders is
$$
\{ 0,1,2,q/2,(q+2)/2,(q+4)/2,q,q+1,q+2,3q/2,3q/2,2q\}\, ,
$$
we conclude that $2q-m_4=q/2+2$; whence $m_4=3q/2-2$. Finally from
(\ref{eq4.14}), $\ell:=\frac{q-4}{3}(q/2+1)+1\not\in H(P)$. On the
other hand, $\ell=m_4+\frac{q-10}{6}m_2\in H(P)$, a contradiction.

(2) $\text{(a)$\Rightarrow$(b)}:$ In virtue of Lemma
\ref{lemma4.11}, we have $g=c_0(q+1,4)$. By
$q+1=\frac{q-2}{3}\cdot 3+3$, Lemma \ref{lemma2.12} shows that
$\dim(2\cD_\cX)=11$ and that $\frac{q-5}{3}\cD_\cX+\cD'$ is the
canonical linear series, where $\cD'$ is a base-point-free
1-dimensional linear series of degree $(q+1)/3$. Let $P\in\cX$ and
assume $j_2(P)>2$ according to Corollary \ref{cor4.11}. If
$P\in\cX(\fq)$, from Lemma \ref{lemma4.13}(1) the result follows.
Otherwise, $P\not\in\cX(\fq)$, and we have two possibilities
according as $q$ is odd or even (Lemma \ref{lemma4.13}(2)).

Case 1: $q$ is odd $j_2(P)=(q-1)/2$, $j_3(P)=(q+1)/2$. A similar
property to (\ref{eq4.14}) holds, namely $\delta+1\not\in H(P)$ for
any $(\frac{q-5}{3}\cD_\cX, P)$-order $\delta$. Hence $2j_2(P)=1=q$
is not in $H(P)$, a contradiction.

Case 2: $q$ is even, $j_2(P)=q/2$, $j_3(P)=(q+2)/2$. From the Case 3
in the proof of part (1), we have $m_1=q/2$. Notice that the degree
$(q+1)/3$ of the above linear series $\cD'$ is coprime to $m_1$.
Then by the well known Riemann's inequality for the genus $g$
applied to $\cD'$ and the linear series corresponding to $m_1$ we
obtain $g\le (q-2)^2/6$, a contradiction.

The implication $\text{(b)$\Rightarrow$(c)}$ is a special case of
\cite[Thm. 2.3]{fgt} while $\text{(c)$\Rightarrow$(a)}$ is trivial.
   \end{proof}

\subsection{The case of $(q-1)(q-3)/8$, $q$ odd}\label{s4.2} The main
result is Theorem \ref{thm4.21} which is analogous to Theorem
\ref{thm4.11}. It states that for $p\ge 5$ and $q$ large enough, the
non-singular $\fq$-model of the curve of equation
   \begin{equation}\label{eq4.21}
y^q+y=x^{(q+1)/4}\, ,\qquad q\equiv 3\pmod{4}\, ,
   \end{equation}
together with the Fermat curve of degree $(q+1)/2$
   \begin{equation}\label{eq4.22}
x^{(q+1)/2}+y^{(q+1)/2}+1=0\, .
   \end{equation}
are the unique $\fq$-maximal curves of genus $g=(q-1)(q-3)/8$
provided that $\dim (\cD_\cX)=5$ holds. The extra-condition on $\dim
(\cD_\cX)$ is assumed since the argument in Lemma \ref{lemma4.11}
only proves that $\dim(\cD_\cX)\in\{4,5\}$. Then
$g=c_0(q+1,5)$, and once again we
take advantage of the hypothesis on the genus by means of Lemma
\ref{lemma2.12}.

The above two curves are in fact not isomorphic even over $\bar\fq$; see
\cite[Remark 4.1]{chkt}. The curve in (\ref{eq4.22}) was characterized in
\cite{chkt} as the unique (up to $\fq$-isomorphism) plane $\fq$-maximal
curve of degree $(q+1)/2$ provided that $q$ is odd and $q\ge 11$.

As $\dim(2\cD_\cX)=14$ by Lemma \ref{lemma2.12}(1), we are able
again to compute the possibilities for the sequence of
$(\cD_\cX,P)$-orders for $P\in\cX$. The proofs of the following two
results will be omited since they are similar to those of Lemmas
\ref{lemma4.12}, \ref{lemma4.13}, and Corollary \ref{cor4.11}. By
Lemma \ref{lemma2.33}(1) $j_1(P)=1$ and either $j_5(P)=q+1$ if
$P\in\cX(\fq)$, or $j_5(P)=q$ otherwise.
   \begin{lemma}\label{lemma4.21} Let $\cX$ be a $\fq$-maximal curve and
$P\in\cX$. Assume that $\dim(\cD_\cX)=5, \dim(2\cD_\cX)=14,$ and
that $q\ge 11.$
   \begin{enumerate}
\item[\rm(1)] If $j_3(P)=3$, then $j_4(P)=4.$
\item[\rm(2)] Let $j_2(P)=2$ but $j_3(P)>3$. If $P\in\cX(\fq),$ then $q$
is odd, $j_3(P)=(q+1)/2,$ and $j_4(P)=(q+3)/2.$  If $P\not\in\cX(\fq),$
then $q$ is even$,$ $j_3(P)=q/2,$ and $j_4(P)=(q+2)/2.$
\item[\rm(3)] Let $P\in\cX(\fq)$ and $j_2(P)>2.$ Assume
$g>(q-2)^2/9$ if $q\equiv 2\pmod{3}$ and $g>(q-3)q/9$ if $q\equiv
0\pmod{3}.$ Then either $q\equiv 3\pmod{4}, j_2(P)=(q+1)/4, j_3(P)=2(q+1)/4,
j_4(P)=3(q+1)/4,$ or $q\equiv 0\pmod{4}, j_2(P)=(q+4)/4, j_3(P)=(2q+4)/4,
j_4(P)=(3q+4)/4.$
\item[\rm(4)] Let $P\not\in\cX(\fq)$ and $j_2(P)>2.$ Then either $q\equiv
1\pmod{4}, j_2(P)=(q+3)/4,$ $j_3(P)=(2q+2)/4, j_4(P)=(3q+1)/4,$ or
$q\equiv 0\pmod{4},$ $j_2(P)=q/4,$ $j_3(P)=2q/4,$ $j_4(P)=3q/4,$ or
$q\equiv 1\pmod{3},$ $j_2(P)=(q-1)/3,$ $j_3(P)=(q+2)/3,$
$j_4(P)=(2q+1)/3,$ or $q\equiv 0\pmod{3}, j_2(P)=q/3, j_3(P)=(q+3)/3,
2q/3.$
   \end{enumerate}
   \end{lemma}
   \begin{corollary}\label{cor4.21} Let $\cX$ be a $\fq$-maximal curve of
genus $g$. Assume that $\dim(\cD_\cX)=5,$ $\dim(2\cD_\cX)=14,$ and
that $q\ge 11$. If $j_3(P)=3$ for every $P\in\cX,$ then $q\equiv 0,
4\pmod{5}$ and $g=(q^2-4q+15)/20.$
   \end{corollary}
   \begin{corollary}\label{cor4.22} Let $\cX$ be a $\fq$-maximal curve of
genus $(q-1)(q-3)/8$ with $q$ odd. Assume $\dim(\cD_\cX)=5$ and $q\ge 11$.
Then:
   \begin{enumerate}
\item[\rm(1)] $\cX$ is $\fq$-isomorphic to the non-singular $\fq$-model of
(\ref{eq4.21}) if and only if there exists $P\in\cX(\fq)$ with $j_2(P)>2;$
\item[\rm(2)] $\cX$ is $\fq$-isomorphic to (\ref{eq4.22}) if and only if
there exists $P\in\cX(\fq)$ with $j_2(P)=2,$ and $j_3(P)>3.$
   \end{enumerate}
   \end{corollary}
   \begin{proof} (1) Let $P$ be the unique point over $x=\infty$. It is
straightforward to check that $m_3(P)=3(q+1)/4$. Hence
$j_2(P)=(q+1)/4$ by Lemma \ref{lemma2.33}(2). Conversely, from
Lemma \ref{lemma4.21}(3) we have $j_4(P)=3(q+1)/4$ and so
$m_1(P)=(q+1)/4$ by Lemma \ref{lemma2.33}(3). Now, the result
follows from \cite[Thm. 2.3]{fgt}.

(2) We have $\cD_\cX=2\cD$, where $\cD$ is the linear series cut
out by lines on $\cX$ (\cite[Thm. 3.5]{chkt}) and hence every
$\fq$-rational inflexion point $P$ (\cite[Lemma 3.6]{chkt})
satisfies both $j_2(P)=2$ and $j_3(P)>3$. Conversely, from Lemmas
\ref{lemma4.21}(2), \ref{lemma2.33}(2) we obtain both
$m_1(P)=(q-1)/2$ and $m_2(P)=(q+1)/2$. Hence the result from
\cite[Thm. 1.1]{chkt}.
   \end{proof}
  \begin{theorem}\label{thm4.21} Let $\cX$ be a $\fq$-maximal curve of
genus $g=(q-1)(q-3)/8$ with $q$ odd$.$ Assume $\dim(\cD_\cX)=5,$ and $p\ge
5.$

    \begin{enumerate}
\item[\rm(1)] If $q\equiv 1\pmod{4}$ and $q\ge 17,$ then $\cX$ is
$\fq$-isomorphic to the Fermat curve (\ref{eq4.22})$.$

\item[\rm(2)] If $q\equiv 3\pmod{4}$ and $q\ge 19,$ then $\cX$ is
$\fq$-isomorphic to either (\ref{eq4.22}) or the
non-singular $\fq$-model of(\ref{eq4.21})$.$
   \end{enumerate}
   \end{theorem}
   \begin{proof} We have already observed that $g=c_0(q+1,5)$
and thus $\dim(2\cD_\cX)=14$. In particular,
by Corollary \ref{cor4.21} there exists $P\in\cX$ with $j_3(P)>3$.

(1) Let $q\equiv 1\pmod{4}$. If $P\in\cX(\fq)$, then Lemma
\ref{lemma4.21}(2)(3) yields $j_2(P)=2$ and the result follows from
Corollary \ref{cor4.22}(2). To show that this is actually the only
possible case, assume on the contrary that $P\not\in\cX(\fq)$. Note
that $\cK:=\frac{q-5}{4}\cD_\cX$ is the canonical linear series by
Lemma \ref{lemma2.12}(2), and hence that $\delta+1\not\in H(P)$ for
any $(\cK,P)$-order $\delta$. Recall that $m_4=q$ by
(\ref{eq2.34}). Now, Lemma \ref{lemma4.21} together with $p\ge 5$
leads to the following two cases.

Case 1: $j_2(P)=(q+3)/4, j_3(P)=(2q+2)/4, j_4(P)=(3q+1)/4.$ Here,
$\{q-m_3, q-m_2, q-m_1\}\subseteq \{1, (q+3)/4, (2q+2)/4,
(3q+1)/4\}$ by Lemma \ref{lemma2.33}(3). Thus $m_1=(2q-2)/4,
m_2=(3q-3), m_3=q-1$. Now, $\delta=(q-9)/4+(3q+1)/4=q-2$ is a
$(\cK,P)$-order and hence $q-1\not\in H(P)$, a contradiction.

Case 2: $q\equiv 1\pmod{3}, j_2(P)=(q-1)/3, j_3(P)=(q+2)/3,
j_4(P)=(2q+1)/3.$ Here, $\delta=3j_2(P)$ is a $(\cK,P)$-order (as
$(q-5)/4\ge 3$) and so $q$ cannot belong to $H(P)$, a
contradiction.

(2) $q\equiv 3\pmod{4}.$ As above, if we show that $P\in\cX(\fq)$, the
result will follow from Corollary \ref{cor4.22}. If $P\not\in\cX(\fq)$,
Lemma \ref{lemma4.21}(2)(4) together with $p\ge 5$ yields
$j_2(P)=(q-1)/3$. Now, Lemma \ref{lemma2.12}(2) implies that
$\delta+1\not\in H(P)$ for every $(\frac{q-7}{4}\cD_\cX,P)$-order
$\delta$. On the other hand, as $(q-7)/4\ge 3$, $3j_2(P)+1=q\in H(P)$, a
contradiction.
   \end{proof}
   \begin{remark}\label{rem4.21} As pointed out in Introduction,
$\fq$-maximal curves of genus $g=\lfloor (q^2-2q+5)/8\rfloor$ do
exist. This genus equals Halphen's number $c_1(4,q+1)$, cf.
(\ref{halphen}). So far, the following examples are known:

   \begin{enumerate}
   \item[\rm(i)] For $q\equiv 0\pmod{4}$, curves of genus $(q^2-2q)/8$
belong to a family of $\fq$-maximal curves constructed by van der
Geer and van der Vlugt, see \cite[Prop. 5.2(ii)]{geer-vl1}, via
fibre products of certain Artin-Schreier $p$-extensions of the
projective line. See also \cite[Thm. 3.3]{g-sti-x}. It seems
plausible that a plane model for such a curve may be obtained from
the proof of \cite[Prop. 1.1]{g-sti}.

   \item[\rm(ii)] For $q\equiv 1\pmod{4}$, curves of genus $(q-1)^2/8$
have been constructed as a quotient of the Hermitian curve $\cH$ by
a subgroup of the automorphism group of $\cH$; see \cite[Prop.
3.3(3)]{ckt2}.

   \item[\rm(iii)] For $q\equiv 3\pmod{4}$, curves of genus $(q^2-2q+5)/8$
have been constructed in a similar way as in (II) above; see \cite[Prop.
3.3(3)(1)]{ckt2} or \cite[Ex. 5.10]{g-sti-x}.
   \end{enumerate}

For the curves mentioned in (ii) and (iii), no plane model seems to
be available in the literature.
   \end{remark}


\begin{thebibliography}{99}

\bibitem{abdon} M. Abd\'on, ``On maximal curves in characteristic two",
Ph.D. dissertation, S\'erie F-121/2000, IMPA, Rio de Janeiro, Brazil,
2000.

\bibitem{at} M. Abd\'on and F. Torres, {\em On maximal curves in
characteristic two}, Manuscripta Math. {\bf 99} (1999), 39-53.

\bibitem{accola} R.D.M. Accola, {\em On Castelnuovo's inequality for
algebraic curves, I}, Trans. Amer. Math. Soc. {\bf 251} (1979), 357--373.

\bibitem{acgh} E. Arbarello, M. Cornalba, P.A. Griffiths, and J. Harris,
``Geometry of Algebraic Curves," Vol. I, Springer-Verlag, New York, 1985.

\bibitem{ballico} E. Ballico, {\em Space curves not contained in low
degree surfaces in positive characteristic}, preprint.

\bibitem{ballico-cossidente} E. Ballico and A. Cossidente, {\em On the
generic hyperplane section of curves in positive characteristic}, J. Pure
and Applied Algebra {\bf 102} (1995) 243-250.

\bibitem{ballico-hefez} E. Ballico and A. Hefez, {\em On the Galois group
associated to a generically etale morphism}, Comm. Algebra {\bf 14}(5)
(1986), 899--909.

\bibitem{castelnuovo} G. Castelnuovo, {\em Ricerche di geometria sulle
curve algebriche}, Atti. R. Acad. Sci. Torino {\bf 24} (1889), 196--223.

\bibitem{chi-ci} L. Chiantini and C. Ciliberto, {\em Towards a Halphen
theory of linera series on curves}, Trans. Amer. Math. Soc. {\bf 351}(6)
(1999), 2197--2212.

\bibitem{chi-ci-ge} L. Chiantini, C. Ciliberto and V. Di Gennaro, {\em The
genus of projective curves}, Duke Math. J. {\bf 70}(2) (1993), 229--245.

\bibitem{chkt} A. Cossidente, J.W.P. Hirschfeld, G. Korchm\'aros and F.
Torres, {\em On plane maximal curves}, Compositio Math. {\bf 121} (2000),
163--181.

\bibitem{ckt1} A. Cossidente, G. Korchm\'{a}ros and F. Torres, {\em
On curves covered by the Hermitian curve}, J. Algebra {\bf 216}
(1999), 56-76.

\bibitem{ckt2} A. Cossidente, G. Korchm\'{a}ros and F. Torres, {\em
Curves of large genus covered by the Hermitian curve}, to appear in Comm.
Algebra.

\bibitem{eisenbud-harris} D. Eisenbud and J. Harris, ``Curves in
projective space'', Les Presses de l'Universit\'e de Montr\'eal,
Montr\'eal, 1982.

\bibitem{fgt} R. Fuhrmann, A. Garcia and F. Torres, {\em On maximal
curves}, J. Number Theory {\bf 67} (1997), 29-51.

\bibitem{ft1} R. Fuhrmann and F. Torres, {\em The genus of curves over
finite fields with many rational points}, Manuscripta Math. {\bf 89}
(1996), 103--106.

\bibitem{ft2} R. Fuhrmann and F. Torres, {\em On Weierstrass points and
optimal curves}, Rend. Circ. Mat. Palermo {\bf 51} (1998), 25--46.

\bibitem{g-q} A. Garcia and L. Quoos, {\em A construction of curves over
finite fields}, preprint May 2000.

\bibitem{g-sti} A. Garcia and H. Stichtenoth, {\em Elementary abelian
$p$-extensions of algebraic function fields}, Manuscripta Math. {\bf 72}
(1991), 67--79.

\bibitem{g-sti1} A. Garcia and H. Stichtenoth, {\em Algebraic functions
fields over finite fields with many places}, IEEE Trans. Inform. Theory
{\bf 41}(6) (1995), 1548--1563.

\bibitem{g-sti-x} A. Garcia, H. Stichtenoth and C.P. Xing, {\em On
subfields of the Hermitian function field}, Compositio Math. {\bf 120}
(2000), 137--170.

\bibitem{gv} A. Garcia and J.F. Voloch, {\em Wronskians and independence
in fields of prime characteristic}, Manuscripta Math. {\bf 59} (1987),
457-469.

\bibitem{geer-vl-1} G. van der Geer and M. van der Vlugt, {\em Weight
distribution for a certain class and maximal curves}, Discr. Math. {\bf
106/107} (1992), 209--218.

\bibitem{geer-vl0} G. van der Geer and M. van der Vlugt, {\em Fibre
products of Artin-Schreier curves and generalized Hamming weights of
codes}, J. Comb. Theory, ser. A {\bf 70} (1995), 337--348.

\bibitem{geer-vl1} G. van der Geer and M. van der Vlugt, {\em Quadratic
forms, generalized Hamming weights of codes and curves with many points},
J. Number Theory {\bf 59} (1996), 20--36.

\bibitem{geer-vl2} G. van der Geer and M. van der Vlugt, How to construct
curves over finite fields with many points, {\em Arithmetic Geometry}
(Cortona 1994) (F. Catanese Ed.), 169--189, Cambridge Univ. Press,
Cambridge, 1997.

\bibitem{geer-vl3} G. van der Geer and M. van der Vlugt, {\em Generalized
Reed-M\"uller codes and curves with many points}, Report W97-22,
Mathematical Institute, University of Leiden, The Netherlands
(alg-geom/9710016).

\bibitem{geer-vl4} G. van der Geer and M. van der Vlugt, {\em Kummer
curves with many points}, preprint math.AG/9909037.

\bibitem{geer-vl} G. van der Geer and M. van der Geer, {\em Tables of
curves with many points}, July 2000,
http:\slash\slash www.wins.uva.nl/\~\ geer.

\bibitem{goppa} V.D. Goppa, ``Geometry and Codes", Math. Appl., Kluwer
Acad. Publ., Dordrecht, 1988.

\bibitem{gruson-peskine} L. Gruson and C. Peskine, Genre des courbes
de l'espace projetif, ``Algebraic Geometry", Proc. Tromso, Norway", Lect.
Notes in Math. Vol. 657, Springer-Verlag, Berlin, 1977.

\bibitem{harris} J. Harris, {\em The genus of space curves}, Math. Ann.
{\bf 249}, (1980), 192--204.

\bibitem{hartshorne} R. Hartshorne, ``Algebraic Geometry", Grad. Texts in
Math., Vol. 52, Springer-Verlag, New York/Berlin, 1977.

\bibitem{hefez-kleiman} A. Hefez and S. Kleiman, Notes on the duality of
projective varieties, ``Geometry Today", 143--183, Birkh\"auser, 1985.

\bibitem{hirschfeld} J.W.P. Hirschfeld, ``Projective Geometries over
Finite Fields", second edition, Oxford University Press, Oxford, 1998.

\bibitem{ihara} Y. Ihara, {\em Some remarks on the number of rational
points of algebraic curves over finite fields}, J. Fac. Sci. Tokyo {\bf
28} (1981), 721--724.

\bibitem{kt1} G. Korchm\'aros and F.Torres, Embedding of a maximal curve
in a Hermitian variety, to appear in Compositio Math.

\bibitem{rathmann} J. Rathmann, {\em The uniform position principle for
curves in characteristic $p$}, Math. Ann. {\bf 276} (1987), 565-579.

\bibitem{r-sti} H.G. R\"uck and  H. Stichtenoth, {\em A characterization
of Hermitian function fields over finite fields}, J. Reine Angew.
Math. {\bf 457} (1994), 185--188.


\bibitem{sti-x} H. Stichtenoth and C. Xing, {\em The genus of maximal
functions fields}, Manuscripta Math. {\bf 86} (1995), 217--224.

\bibitem{sv} K.O. St\"{o}hr and J.F. Voloch, {\em Weierstrass points and
curves over finite fields}, Proc. London Math. Soc. {\bf 52} (1986),
1-19.


\end{thebibliography}
\end{document}